\documentclass[11pt]{article}
\usepackage{amssymb,latexsym}
\input epsf  

\title{On the M\"obius function of a lower Eulerian Cohen-Macaulay poset}

\author{Christos~A.~Athanasiadis\\
Department of Mathematics\\
University of Athens\\
Athens 15784, Hellas (Greece)\\
\small Email: \texttt{caath@math.uoa.gr} }

\date{\small May 16, 2011; revised, July 4, 2011}

  \def\kk{{\mathbf k}}
  \def\aA{{\mathcal A}}
  \def\bB{{\mathcal B}}

  \def\mM{{\mathcal M}}
   
  \def\cost{{\rm cost}}
  \def\inte{{\rm int}}
  
  \def\link{{\rm lk}} 
  \def\star{{\rm st}}

  \def\sm{\smallsetminus} 
  
  \newcommand{\qed}{$\hfill \Box$}

\begin{document}
\maketitle

\newtheorem{theorem}{Theorem}[section]
\newtheorem{proposition}[theorem]{Proposition}
\newtheorem{corollary}[theorem]{Corollary}
\newtheorem{definition}[theorem]{Definition}
\newtheorem{remark}[theorem]{Remark}
\newtheorem{lemma}[theorem]{Lemma}
\newtheorem{example}[theorem]{Example}
\newtheorem{examples}[theorem]{Examples}
\newtheorem{conjecture}[theorem]{Conjecture}
\newtheorem{question}[theorem]{Question}

\begin{abstract}
A certain inequality is shown to hold for the values of the M\"obius 
function of the poset obtained by attaching a maximum element to a lower 
Eulerian Cohen-Macaulay poset. In two important special cases, this 
inequality provides partial results supporting Stanley's nonnegativity 
conjecture for the toric $h$-vector of a lower Eulerian Cohen-Macaulay 
meet-semilattice and Adin's nonnegativity conjecture for the cubical 
$h$-vector of a Cohen-Macaulay cubical complex.     

\end{abstract}

  \section{Introduction}
  \label{sec:intro}

  Let $P$ be a finite poset which has a minimum element $\hat{0}$ (for background and 
  any undefined terminology on partially ordered sets we refer the reader to \cite[Chapter 
  3]{StaEC1} and Section \ref{sec:pre}). Such a poset is called \emph{lower Eulerian} 
  \cite[Section 4]{Sta87} if every interval $[x, y]$ in $P$ is graded and satisfies $\mu_P 
  (x, y) = (-1)^{\rho(x, y)}$, where $\mu_P$ is the M\"obius function of $P$ and $\rho(x, 
  y)$ is the common length of all maximal chains of $[x, y]$. Examples of lower Eulerian 
  posets are the face posets of finite regular cell complexes; see \cite[Section 
  12.4]{Bj95} \cite[Section 3.8]{StaEC1}. The main result of this paper concerns the 
  M\"obius function of a poset obtained by attaching a maximum element to a lower Eulerian 
  Cohen-Macaulay poset (here, the Cohen-Macaulay property is defined with respect to an 
  arbitrary field).  

  \begin{theorem} \label{thm:main}
    Let $P$ be a lower Eulerian Cohen-Macaulay poset, with minimum element $\hat{0}$ and 
    set of atoms $\aA(P)$. Let $\hat{P} = P \cup \{\hat{1}\}$ be the poset obtained from 
    $P$ by attaching a maximum element $\hat{1}$ and let $\mu_{\hat{P}}$ denote the M\"obius 
    function of $\hat{P}$. Then we have
      \begin{equation} \label{eq:main}
        \sum_{x \in \aA(P)} \, |\mu_{\hat{P}} (x, \hat{1})| \ \ge \ \alpha_P \, 
        |\mu_{\hat{P}} (\hat{0}, \hat{1})|,
      \end{equation} 
    where $\alpha_P$ is the minimum cardinality of the sets $\{x \in \aA(P): x \le_P 
    y\}$ of atoms of $P$ in the interval $[\hat{0}, y]$, when $y$ runs through all 
    maximal elements of $P$.  
  \end{theorem}

  To explain the motivation behind Theorem \ref{thm:main} and discuss some of its 
  consequences, we consider two important classes of lower Eulerian posets. We recall 
  that a finite poset $P$ having a minimum element $\hat{0}$ is called \emph{simplicial} 
  \cite[p.~135]{StaEC1} (respectively, \emph{cubical} \cite{BBC97}) if for every $x \in 
  P$, the interval $[\hat{0}, x]$ in $P$ is isomorphic to a Boolean lattice 
  (respectively, to the face poset of a cube). The $h$-vector $(h_0 (P), h_1 (P),\dots,h_d 
  (P))$ is a fundamental enumerative invariant of a simplicial poset $P$ (see, for 
  instance, \cite[p.~113]{Sta87}) defined by the formula 
    \begin{equation} \label{eq:h}
      \sum_{i=0}^d \, h_i (P) q^i \ = \ \sum_{i=0}^d \, f_{i-1} (P) \, q^i (1-q)^{d-i},
    \end{equation}
  where $d$ is the largest length of a chain in $P$ and $f_{i-1} (P)$ is the number of 
  elements $x \in P$ such that the interval $[\hat{0}, x]$ is isomorphic to a Boolean 
  lattice of rank $i$. For a graded cubical poset $P$ of rank $d$, Adin \cite{Ad96}  
  introduced the cubical $h$-vector $(h^{(c)}_0 (P), h^{(c)}_1 (P),\dots,h^{(c)}_d (P))$ 
  of $P$ as a cubical analogue of the $h$-vector of a simplicial poset. It can be defined 
  by the formula
    $$ \sum_{i=0}^d h^{(c)}_i (P) q^i \ = \ \frac{1}{1+q} \, \left( 2^{d-1} \, + \, 
      q \, \sum_{i=0}^{d-1} \, f_i (P) \, (2q)^i (1-q)^{d-i-1} \, + \, 
      (-2)^{d-1} \, \widetilde{\chi} (P) \, q^{d+1} \right), $$
  where $f_i (P)$ denotes the number of elements of $P$ of rank $i+1$, for $0 \le i \le 
  d-1$, and $\widetilde{\chi} (P) = - 1 + \sum_{i=0}^{d-1} (-1)^i f_i (P)$. 

  \medskip
  One of Stanley's early results in $f$-vector theory (see \cite[Corollary 4.3]{Sta75}) 
  states that the $h$-vector of $P$ has nonnegative entries if $P$ is the face poset of a 
  Cohen-Macaulay simplicial complex. This statement was extended to all Cohen-Macaulay 
  simplicial posets in \cite{Sta91} and, in fact, the $h$-vectors of Cohen-Macaulay 
  simplicial complexes and posets can both be completely characterized \cite{Sta77, Sta91} 
  (see also Sections II.3 and III.6 in \cite{StaCCA}). On the contrary, far less is known 
  about $f$-vectors of cubical complexes and posets. Adin \cite[Question 1]{Ad96} (see 
  also \cite[Problem 8 (a)]{Sta00}) raised the question whether the cubical $h$-vector of 
  $P$ has nonnegative entries if $P$ is the face poset of a Cohen-Macaulay cubical complex 
  and gave an affirmative answer in the special case of shellable cubical complexes 
  \cite[Theorem 5 (iii)]{Ad96}. For the importance of Adin's question we refer the reader 
  to \cite{Het96}, where an affirmative answer is given for cubical barycentric 
  subdivisions (implicitly defined in Section \ref{subsec:inte}) of simplicial complexes. 

  One can easily check directly that $h_d (P) \ge 0$ (respectively, $h^{(c)}_d (P) \ge 0$) 
  holds for every simplicial (respectively, cubical) Cohen-Macaulay poset $P$ of rank $d$. 
  Using the notation of Theorem \ref{thm:main}, we claim (see Section \ref{sec:app}) that  
    \begin{equation} \label{eq:mains}
      h_{d-1} (P) \ = \ (-1)^d \sum_{x \in \aA(P)} \mu_{\hat{P}} (x, \hat{1}) \, - \, d \, 
      (-1)^{d+1} \mu_{\hat{P}} (\hat{0}, \hat{1})
    \end{equation} 
  and

    \begin{equation} \label{eq:mainc}
      h^{(c)}_{d-1} (P) \ = \ (-1)^d \sum_{x \in \aA(P)} \mu_{\hat{P}} (x, \hat{1}) \, - \, 
      2^{d-1} \, (-1)^{d+1} \mu_{\hat{P}} (\hat{0}, \hat{1})
    \end{equation} 
  holds for every simplicial (respectively, cubical) graded poset $P$ of rank $d$. Thus,
  Theorem \ref{thm:main} applies to both of these situations and gives the following partial
  information on the nonnegativity of simplicial and cubical $h$-vectors. Since the proof 
  of Theorem \ref{thm:main} uses tools from topological combinatorics, it specializes to a 
  new proof that $h_{d-1} (P) \ge 0$ holds for Cohen-Macaulay simplicial posets $P$ of rank 
  $d$. The case of cubical posets yields a new result as follows.
    \begin{corollary} \label{cor:cubical}
      For every Cohen-Macaulay cubical poset $P$ of rank $d$ we have $h^{(c)}_{d-1} (P) 
      \ge 0$. 
    \end{corollary}
  Motivated by results on the intersection cohomology of toric varieties, Stanley 
  \cite[Section 4]{Sta87} defined the (generalized, or toric) $h$-vector $(h_0 (P), h_1 
  (P),\dots,h_d (P))$ for an arbitrary lower Eulerian poset $P$ (where again, $d$ is the 
  largest length of a chain in $P$). For simplicial posets, this $h$-vector reduces to 
  the one defined by (\ref{eq:h}). Stanley \cite[Conjecture 4.2 (b)]{Sta87} conjectured 
  that the generalized $h$-vector has nonnegative entries for every lower Eulerian 
  Cohen-Macaulay meet-semilattice. The following partial result is also a consequence of 
  Theorem \ref{thm:main} (see Section \ref{sec:toric}).   
    \begin{corollary} \label{cor:toric}
      For every lower Eulerian Cohen-Macaulay meet-semilattice $P$ of rank $d$ we have 
      $h_{d-1} (P) \ge 0$. 
    \end{corollary}
  This paper is organized as follows. Section \ref{sec:pre} reviews basic definitions and
  background on the enumerative and topological combinatorics of partially ordered sets and 
  establishes some preliminary results on (lower Eulerian) Cohen-Macaulay posets. Section 
  \ref{sec:buch} proves a certain statement (Corollary \ref{cor:surjectivity}) on upper
  truncations of lower Eulerian Cohen-Macaulay posets which will be essential in the 
  proof of Theorem \ref{thm:main}. This statement follows from more general statements on 
  rank-selections of Cohen-Macaulay (or Buchsbaum) posets and balanced simplicial 
  complexes (Theorems \ref{thm:2CMposets} and \ref{thm:2CM}), essentially established
  by Browder and Klee in \cite{BK}. Theorem \ref{thm:main} is proved in Section 
  \ref{sec:proof}. The applications to cubical and toric $h$-vectors are discussed in 
  Section \ref{sec:app}.

  \section{Preliminaries}
  \label{sec:pre}

  This section reviews basic background on the enumerative and topological combinatorics
  of simplicial complexes and partially ordered sets (posets), fixes notation and 
  establishes some preliminary results which will be useful in the sequel. For more 
  information on these topics we refer the reader to \cite[Chapter 3]{StaEC1} and 
  \cite{Bj95}. Basic background on algebraic topology can be found in \cite{Mu84}. 

  \subsection{Simplicial complexes}
  \label{subsec:complexes}

  Given a finite set $E$, an (abstract) \emph{simplicial complex} on the ground set $E$ 
  is a collection $\Delta$ of subsets of $E$ such that $\sigma \subseteq \tau \in \Delta$ 
  implies $\sigma \in \Delta$. The elements of $\Delta$ are called \emph{faces}. The 
  dimension of a face $\sigma$ is defined as one less than the cardinality of $\sigma$. 
  A \emph{facet} of $\Delta$ is a face which is maximal with respect to inclusion. The 
  \emph{reduced Euler characteristic} of $\Delta$ is defined as  
    \begin{equation} \label{eq:chiDelta}
      \widetilde{\chi} (\Delta) \ = \ \sum_{i=0}^d \, (-1)^{i-1} f_{i-1} (\Delta),  
    \end{equation} 
  where $d-1 = \dim \Delta$ is the dimension (maximum dimension of a face) of $\Delta$ 
  and $f_{i-1} (\Delta)$ denotes the number of faces of $\Delta$ of dimension $i-1$, for 
  $0 \le i \le d$. By the Euler-Poincar\'e formula we have 
    \begin{equation} \label{eq:EP}
      \widetilde{\chi} (\Delta) \ = \ \sum_{i=0}^d \, (-1)^{i-1} \, \dim_\kk 
      \widetilde{H}_{i-1} (\Delta; \kk),
    \end{equation} 
  where $\widetilde{H}_* (\Delta; \kk)$ stands for the reduced simplicial homology of
  $\Delta$ over the field $\kk$. The \emph{link} of a face $\sigma \in \Delta$ is defined 
  as the subcomplex $\link_\Delta (\sigma) = \{\tau \sm \sigma: \tau \in \Delta, \, 
  \sigma \subseteq \tau\}$ of $\Delta$. The simplicial complex $\Delta$ is called 
  \emph{Cohen-Macaulay} over $\kk$ if $\widetilde{H}_i \, (\link_\Delta (\sigma); \kk) 
  = 0$ for every $\sigma \in \Delta$ (including $\sigma = \varnothing$) and all $i < 
  \dim \link_\Delta (\sigma)$. Such a complex is \emph{pure}, meaning that every facet 
  of $\Delta$ has dimension equal to $\dim \Delta$. 
  All topological properties of $\Delta$ we mention in the sequel will refer to those 
  of the geometric realization \cite[Section 9]{Bj95} of $\Delta$, uniquely defined up 
  to homeomorphism.

  \subsection{The M\"obius function and the order complex}
  \label{subsec:mobius}

  A poset $P$ is called \emph{locally finite} if every closed interval $[x, y]$ in $P$ is 
  finite. The \emph{M\"obius function} $\mu_P$ of such a poset is defined on pairs $(x, y)$
  of elements of $P$ satisfying $x \le_P y$ by the recursive formula 
    \begin{equation} \label{eq:defmobius}
      \mu_P (x, y) \ = \ \cases{ 1, & if \ $x = y$, \cr
      {\displaystyle \ - \sum_{x \, \le_P \, z \, <_P \, y} \mu_P (x, z)}, & 
      if \ $x <_P y$. } 
    \end{equation} 
  Given a finite poset $Q$, we denote by $\Delta (Q)$ the simplicial complex on the ground 
  set $Q$ whose faces are chains (totally ordered subsets) of $Q$, known as the \emph{order 
  complex} of $Q$; see \cite[p.~120]{StaEC1} \cite[Section 9]{Bj95}. The following 
  proposition gives a fundamental interpretation of the M\"obius function of a locally 
  finite poset.

  \begin{proposition} {\rm (\cite[Proposition 3.8.6]{StaEC1})} \label{prop:hall}
    For every locally finite poset $P$ and all $x, y \in P$ with $x <_P y$ we have 
    $\mu_P (x, y) = \widetilde{\chi} (\Delta(x, y))$, where $\Delta(x, y)$ denotes the 
    order complex of the open interval $(x, y)$ in $P$. 
  \end{proposition}

  \subsection{Lower Eulerian posets}
  \label{subsec:lower}

  A locally finite poset $P$ is called \emph{locally graded} if for every closed interval 
  $[x, y]$ in $P$ there exists a nonnegative integer $\rho(x, y)$, called the rank of $[x, 
  y]$, such that every maximal chain in $[x, y]$ has length equal to $\rho(x, y)$. A locally 
  graded poset which has a minimum element will be called \emph{lower graded}. For such a 
  poset $P$ with minimum element $\hat{0}$, we will refer to the rank of 
  $[\hat{0}, x]$ simply as the rank of $x$ and will denote it by $\rho(x)$. A locally 
  graded poset $P$ is called \emph{locally Eulerian} if $\mu_P (x, y) = (-1)^{\rho (x, 
  y)}$ holds for all $x, y \in P$ with $x \le_P y$. A \emph{lower Eulerian} (respectively, 
  \emph{Eulerian}) poset is a locally Eulerian poset which has a minimum (respectively, 
  a minimum and a maximum) element. 

  Given a finite poset $P$ with a minimum element $\hat{0}$, we will denote by $\bar{P}$ 
  (respectively, by $\hat{P}$) the poset which is obtained from $P$ by removing $\hat{0}$ 
  (respectively, by attaching a maximum element $\hat{1}$). We set 
    \begin{equation} \label{eq:chi}
      \widetilde{\chi} (P) \ = \ \mu_{\hat{P}} (\hat{0}, \hat{1}) 
                           \ = \ \widetilde{\chi} (\Delta (\bar{P})), 
    \end{equation} 
  where the second equality is due to Proposition \ref{prop:hall}. For a finite, lower 
  graded poset $P$ we also set 
    \begin{equation} \label{eq:psi}
      \psi (P) \ = \ \sum_{i=0}^d \, (-1)^{i-1} f_{i-1} (P), 
    \end{equation} 
  where $d$ is the maximum rank of an element of $P$ and $f_{i-1} (P)$ denotes the number 
  of elements of $P$ of rank $i$, for $0 \le i \le d$.
    \begin{lemma} \label{lem:chi}
      We have $\widetilde{\chi} (P) = \psi (P)$ for every lower Eulerian poset $P$.
    \end{lemma}

    \noindent
    \emph{Proof.}
      Using the defining recurrence (\ref{eq:defmobius}) for the M\"obius function and 
      the fact that each closed interval in $P$ is Eulerian, we find that
        $$ \mu_{\hat{P}} (\hat{0}, \hat{1}) \ = \ - \ \sum_{x \in P} \ \mu_P (\hat{0}, x) 
           \ = \ \sum_{x \in P} \ (-1)^{\rho (x)-1} $$
      and the proof follows.
    \qed

  \medskip
  The \emph{rank} of an Eulerian poset is defined as the rank of its maximum element. 
  The following elementary lemma will be used in Section \ref{sec:proof}.
  
  \begin{lemma} \label{lem:mobiusminusatoms}
    Let $P$ be an Eulerian poset of rank $d \ge 2$ with minimum element $\hat{0}$ and 
    maximum element $\hat{1}$ and let $Q$ be the poset which is obtained from $P$ by 
    removing all its atoms. Then we have
      $$ \mu_Q (\hat{0}, \hat{1}) \ = \ (-1)^{d-1} (f_0 (P) - 1), $$ 
  where $f_0 (P)$ is the number of atoms of $P$.
  \end{lemma}

  \noindent
  \emph{Proof.}
    This statement is a special case of \cite[Proposition 3.14.5]{StaEC1}.
  \qed

  \subsection{The poset of intervals}
  \label{subsec:inte}

  We denote by $\inte(P)$ the set of (nonempty) closed intervals $[x, y]$ of a locally
  finite poset $P$, partially ordered by inclusion. For the enumerative, order-theoretic
  and topological properties of $\inte(P)$ we refer the reader to Exercises 7 and 58 of 
  \cite[Chapter 3]{StaEC1} and to \cite[Section 2]{BBC97} \cite[Section 4]{Wa88}. The 
  following proposition (in the form of \cite[Theorem 4.1]{Wa88}) also appears as 
  Theorem 2.3 in \cite{BBC97}.  

  \begin{proposition} {\rm (\cite[Theorem 6.1 (a)]{Wa88})} \label{prop:inte}
    For every finite poset $P$, the order complex $\Delta(\inte(P))$ is homeomorphic
    to $\Delta(P)$.
  \end{proposition}

  \begin{proposition} \label{prop:inteLE}
    For every locally Eulerian poset $P$, the poset which is obtained from $\inte(P)$ by
    attaching a minimum element is lower Eulerian.
  \end{proposition}

  \noindent
  \emph{Proof.}
    Let $Q$ be the poset in question. We denote by $\hat{0}$ the minimum element of $Q$ 
    and choose any elements $x, y \in Q$ with $x <_Q y$. Assume first that $x \ne \hat{0}$. 
    Then we may write $x = [b, c] \in \inte(P)$ and $y = [a, d] \in \inte(P)$ with $a \le_P 
    b \le_P c \le_P d$. Clearly, the interval $[x, y]$ in $Q$ is isomorphic to the direct 
    product of the intervals $[a, b]$ and $[c, d]$ in $P$ and hence we have $\mu_Q (x, y) 
    = \mu_P (a, b) \mu_P (c, d)$ by \cite[Proposition 3.8.2]{StaEC1}. Suppose now that $x 
    = \hat{0}$ and let $y = [a, d] \in \inte(P)$, as before. Then $y$ (as a subposet of $P$) 
    is an Eulerian poset and the interval $[x, y]$ in $Q$ is the poset obtained from 
    $\inte(y)$ by attaching a minimum element $\hat{0}$. By the result of \cite[Excercise 
    3.58 (b)]{StaEC1} we have $\mu_Q (x, y) = \mu_Q (\hat{0}, y) = - \mu_P (a, d)$. The 
    previous observations and the fact that each closed interval in $P$ is Eulerian imply 
    that $Q$ is locally Eulerian as well.  
  \qed

  \subsection{Cohen-Macaulay posets}
  \label{subsec:CM}

  A finite poset $P$ is \emph{graded} of rank $d$ (respectively, Cohen-Macaulay over the 
  field $\kk$) if the order complex $\Delta(P)$ is pure of dimension $d$ (respectively, 
  Cohen-Macaulay over $\kk$). 
  For the remainder of this section we assume that $P$ has a minimum element $\hat{0}$. 
  Then $P$ is Cohen-Macaulay over $\kk$ if and only if $\bar{P}$ is as well. Moreover, 
  in that case $\bar{P}$ and $P$ are graded of rank $d-1$ and $d$, respectively, and by 
  (\ref{eq:EP}) and (\ref{eq:chi}) we have
    \begin{equation} \label{eq:chiCM}
      \widetilde{\chi} (P) \ = \ (-1)^{d-1} \, \dim_\kk \widetilde{H}_{d-1} (\Delta 
      (\bar{P}); \kk).
    \end{equation} 
  Suppose now that $P$ is graded of rank $d \ge 2$ and let $Q$ be the poset which is 
  obtained from $P$ by removing the set $\mM(P)$ of maximal elements of $P$ (to follow 
  the proofs in the present and the following two sections, it may be helpful for the 
  reader to keep in mind the special case in which $P$ is the face poset of a regular 
  cell complex). Given $y \in \mM(P)$, the order complex $\Delta (\hat{0}, y)$ of the 
  open interval $(\hat{0}, y)$ of $P$ is a subcomplex of $\Delta(\bar{Q})$ and hence 
  there is a map $\widetilde{H}_{d-2} (\Delta (\hat{0}, y); \kk) \to \widetilde{H}_{d-2} 
  (\Delta(\bar{Q}); \kk)$, induced by inclusion. We denote by $\Omega(y)$ the image of 
  this map. Equations (\ref{eq:chi}) and (\ref{eq:chiCM}) imply that if $P$ is lower 
  Eulerian and Cohen-Macaulay over $\kk$, then $\Omega(y)$ is a one-dimensional 
  $\kk$-vector space for every $y \in \mM(P)$. We will then denote by $\omega(y)$ any 
  basis (nonzero) element of $\Omega(y)$.  
 
  \begin{lemma} \label{lem:Qhomology}
    Suppose that $\bar{P}$ is Cohen-Macaulay over $\kk$ of rank $d-1 \ge 1$ and let 
    $\bar{Q}$ be the poset which is obtained from $\bar{P}$ by removing all maximal 
    elements. Then the $\kk$-vector space $\widetilde{H}_{d-2} (\Delta(\bar{Q}); \kk)$ 
    is equal to the sum of its subspaces $\Omega(y)$ for $y \in \mM(\bar{P})$. 
  \end{lemma}   
  
  \noindent
  \emph{Proof.}
    Let $c$ be a $(d-2)$-cycle in the chain complex of $\Delta(\bar{Q})$ over $\kk$. We 
    will show that $c$ can be written as a sum of $(d-2)$-cycles in the chain complexes 
    of $\Delta (\hat{0}, y)$ over $\kk$, for $y \in \mM(\bar{P})$. We observe that
    $c$ is also a $(d-2)$-cycle in the chain complex of $\Delta (\bar{P})$ over $\kk$. 
    Since $\Delta (\bar{P})$ is Cohen-Macaulay over $\kk$ of dimension $d-1$, we have 
    $\widetilde{H}_{d-2} (\Delta(\bar{P}); \kk) = 0$ and hence $c = \partial_{d-1} 
    (\tilde{c})$ for some $(d-1)$-chain $\tilde{c}$ in the chain complex of $\Delta 
    (\bar{P})$ over $\kk$, where $\partial_*$ denotes the boundary map of this complex. 
    Clearly $\tilde{c}$ can be written (uniquely) as the sum of $(d-1)$-chains 
    $\tilde{c}_y$ for $y \in \mM(\bar{P})$, where $\tilde{c}_y$ is a $(d-1)$-chain in 
    the chain complex of $\Delta (\hat{0}, y]$ over $\kk$. Then $c$ is equal to the sum 
    of the boundaries $\partial_{d-1} (\tilde{c}_y)$. Since $c$ is supported on $\Delta
    (\bar{Q})$, each chain $\partial_{d-1} (\tilde{c}_y)$ is supported on $\Delta (\hat{0}, 
    y)$. Since $\partial_{d-2} \, \partial_{d-1} = 0$, this implies that $\partial_{d-1} 
    (\tilde{c}_y)$ is a $(d-2)$-cycle in the chain complex of $\Delta (\hat{0}, y)$ 
    over $\kk$. This proves the desired statement for $c$.
  \qed

  \begin{corollary} \label{cor:Qhomology}
    Let $P$ be a lower Eulerian Cohen-Macaulay poset of rank $d \ge 2$ and let $Q$ denote 
    the poset which is obtained from $P$ by removing all maximal elements. Then the 
    $\kk$-vector space $\widetilde{H}_{d-2} (\Delta(\bar{Q}); \kk)$ is spanned by the 
    classes $\omega(y)$ for $y \in \mM(P)$. 
  \end{corollary}
  
  \noindent
  \emph{Proof.}
    This statement is a special case of Lemma \ref{lem:Qhomology}.
  \qed

  \section{Doubly Cohen-Macaulay and Buchsbaum* posets}
  \label{sec:buch}

  This section proves a certain property (Corollary \ref{cor:surjectivity}) of the first 
  upper truncation $Q$ of a lower Eulerian Cohen-Macaulay poset which will be essential 
  in the proof of Theorem \ref{thm:main} in Section \ref{sec:proof}. This property follows 
  from the Buchsbaum* condition, recently introduced in \cite{AW}, for the order complex  
  $\Delta (\bar{Q})$, so we begin by recalling the relevant definitions. Throughout this 
  section, we write $H_i (\Delta, \Gamma; \kk)$ for the relative simplicial homology of 
  the pair of simplicial complexes $(\Delta, \Gamma)$, where $\Gamma$ is a subcomplex of 
  $\Delta$. 

  A simplicial complex $\Delta$ is \emph{Buchsbaum} over the field $\kk$ if $\Delta$ 
  is pure and $\widetilde{H}_i \, (\link_\Delta (\sigma); \kk) = 0$ for every nonempty face 
  $\sigma \in \Delta$ and all $i < \dim \link_\Delta (\sigma)$. Recall that the 
  \emph{contrastar} of a face $\sigma \in \Delta$ is defined as the subcomplex $\cost_\Delta 
  (\sigma) = \{\tau \in \Delta: \sigma \not\subseteq \tau\}$ of $\Delta$. 

  \begin{proposition} {\rm (cf. \cite[Theorem 8.1]{StaCCA})} \label{prop:buchchar}
    The following conditions on a simplicial complex $\Delta$ are equivalent:
    \begin{itemize}
      \item[(i)] $\Delta$ is Buchsbaum over $\kk$.
      \item[(ii)] $\Delta$ is pure and $\link_\Delta(\sigma)$ is Cohen-Macaulay 
                    over $\kk$ for every $\sigma \in \Delta \sm \{ \varnothing \}$. 
      \item[(iii)] $H_i (\Delta, \cost_\Delta(\sigma); \kk) = 0$ holds for every $\sigma 
                   \in \Delta \sm \{ \varnothing \}$ and all $i < \dim \Delta$.
    \end{itemize}
  \end{proposition}

  The following definition of a Buchsbaum* simplicial complex is equivalent to 
  \cite[Definition 1.2]{AW} (see Propositions 2.3 and 2.8 in \cite{AW}) and will be
  convenient for the purposes of this section.

  \begin{definition} \label{def:Buch*}
    Let $\Delta$ be a $(d-1)$-dimensional simplicial complex which is Buchsbaum over 
    $\kk$. The complex $\Delta$ is called \emph{Buchsbaum* over $\kk$} if the canonical 
    map
      $$ \rho_*: \widetilde{H}_{d-1} (\Delta; \kk) \rightarrow H_{d-1}(\Delta, 
                  \cost_\Delta(\sigma); \kk), $$
    induced by inclusion, is surjective for every $\sigma \in \Delta \sm \{\varnothing\}$. 
  \end{definition}

  Assume now that $\Delta$ is Cohen-Macaulay over $\kk$. Then $\Delta$ is called 
  \emph{doubly Cohen-Macaulay} over $\kk$ \cite[p.~71]{StaCCA} if for every vertex $v$ of 
  $\Delta$, the complex $\Delta \sm v = \{\tau \in \Delta: v \not\in \tau\}$ (obtained 
  from $\Delta$ by removing all faces which contain $v$) is Cohen-Macaulay over $\kk$ of 
  the same dimension as $\Delta$. Given that $\Delta$ is Cohen-Macaulay over $\kk$, it
  was shown in \cite[Proposition 2.5]{AW} that $\Delta$ is doubly Cohen-Macaulay over 
  $\kk$ if and only if $\Delta$ is Buchsbaum* over $\kk$. A Cohen-Macaulay complex 
  $\Delta$ over $\kk$ is called \emph{Gorenstein*} over $\kk$ \cite[Section II.5]{StaCCA} 
  if for every $\sigma \in \Delta$ and for $i = \dim \link_\Delta (\sigma)$ we have 
  $\widetilde{H}_i \, (\link_\Delta (\sigma); \kk) \cong \kk$. Every such complex is 
  doubly Cohen-Macaulay over $\kk$ (see, for instance, Theorem II.5.1 and Proposition 
  III.3.7 in \cite{StaCCA}). A finite poset $Q$ is \emph{Buchsbaum} (respectively, 
  \emph{Buchsbaum*} or \emph{doubly Cohen-Macaulay} or \emph{Gorenstein*}) over $\kk$ 
  if $\Delta(Q)$ is a Buchsbaum (respectively, Buchsbaum* or doubly Cohen-Macaulay or 
  Gorenstein*) simplicial complex over $\kk$. 

  Part (i) of the next theorem follows from 
  \cite[Corollary 2.7]{Fl05} in the special case in which $P$ is the face poset of a 
  regular cell complex with the intersection property and from \cite[Theorem 4.5]{Ya} 
  in the special case in which $P$ is a simplicial poset.

  \begin{theorem} \label{thm:2CMposets}
    Let $P$ be a graded poset of rank $d \ge 2$ with a minimum element $\hat{0}$ and let 
    $Q$ be the poset which is obtained from $P$ by removing all maximal elements. 
      \begin{itemize}
        \item[(i)] If $P$ is Cohen-Macaulay over $\kk$ and the interval $(\hat{0}, y)$ 
        in $P$ is doubly Cohen-Macaulay over $\kk$ for every maximal element $y$ of $P$, 
        then $\bar{Q}$ is doubly Cohen-Macaulay over $\kk$.
        \item[(ii)] If $\bar{P}$ is Buchsbaum over $\kk$ and the interval $(\hat{0}, y)$ 
        in $P$ is doubly Cohen-Macaulay over $\kk$ for every maximal element $y$ of $P$, 
        then $\bar{Q}$ is Buchsbaum* over $\kk$.
      \end{itemize}
  \end{theorem}

  Before we comment on the proof of Theorem \ref{thm:2CMposets}, let us deduce the 
  statement which will be needed in Section \ref{sec:proof}. Let $\Gamma$ be a 
  simplicial complex of dimension $d-2$. Recall that the \emph{closed star} of a vertex 
  $v$ in $\Gamma$ is the subcomplex of $\Gamma$ defined as $\overline{\star}_\Gamma (v) 
  = \{\sigma \in \Gamma: \sigma \cup \{v\} \in \Gamma\}$. Via the isomorphisms
  $H_i \, (\Gamma, \cost_\Gamma (v); \kk) \cong H_i \, (\overline{\star}_\Gamma (v), 
  \link_\Gamma (v); \kk) \cong \widetilde{H}_{i-1} (\link_\Gamma (v); \kk)$, the 
  canonical map $\rho_*: \widetilde{H}_{d-2} (\Gamma; k) \rightarrow H_{d-2}(\Gamma, 
  \cost_\Gamma (v); \kk)$, considered in Definition \ref{def:Buch*}, induces a map 
    $$ \rho_v: \widetilde{H}_{d-2} (\Gamma; \kk) \rightarrow 
               \widetilde{H}_{d-3} (\link_\Gamma (v); \kk). $$
  One can check that $\rho_v$ is induced by a chain map from the (augmented) chain 
  complex of $\Gamma$ over $\kk$ to that of $\link_\Gamma (v)$ which sends a face 
  $\sigma \in \Gamma$ to $\sigma \sm \{v\}$ (with appropriate sign), if $v \in \sigma$, 
  and to zero otherwise.
  Suppose now that $\Gamma = \Delta (\bar{Q})$, where $Q$ is a graded poset having a 
  minimum element and rank $d-1$, and let $x \in Q$ be an atom. Since $\link_\Gamma (x) 
  = \Delta (Q_{> x})$, where $Q_{> x} = \{y \in Q: x <_Q y\}$ is considered as a 
  subposet of $Q$, we get a map
    $$ \rho_x: \widetilde{H}_{d-2} (\Delta(\bar{Q}); \kk) \rightarrow 
               \widetilde{H}_{d-3} (\Delta(Q_{> x}); \kk). $$

  \begin{corollary} \label{cor:surjectivity}
    Let $P$ be a lower Eulerian Cohen-Macaulay poset (over $\kk$) of rank $d \ge 2$ 
    and let $\bar{Q}$ denote the poset which is obtained from $P$ by removing the 
    minimum and all maximal elements. Then the map $\rho_x: \widetilde{H}_{d-2} 
    (\Delta(\bar{Q}); \kk) \rightarrow \widetilde{H}_{d-3} (\Delta(Q_{> x}); \kk)$ 
    is surjective for every minimal element $x$ of $\bar{Q}$. 
  \end{corollary}   

  \noindent
  \emph{Proof.}
    Our assumptions imply that the interval $(\hat{0}, y)$ in $P$ is Gorenstein*, and
    hence doubly Cohen-Macaulay, over $\kk$ for every maximal element $y$ of $P$. Thus,
    Theorem \ref{thm:2CMposets} implies that $\Delta(\bar{Q})$ is Buchsbaum* over $\kk$. 
    In view of Definition \ref{def:Buch*}, the proof follows from the discussion 
    preceding the statement of the corollary. 
  \qed

  \medskip
  Theorem \ref{thm:2CMposets} will be deduced from the following more general statement
  on simplicial complexes. The proof of the latter will follow from that of \cite[Theorem 
  3.1]{BK}. For a subset $U$ of the set of vertices of $\Delta$ we write $\Delta \sm U$ 
  for the subcomplex $\{\sigma \in \Delta: \sigma \cap U = \varnothing\}$ of $\Delta$, 
  obtained from $\Delta$ by removing all faces which intersect $U$.

  \begin{theorem} {\rm (cf. \cite[Theorem 3.1]{BK})} \label{thm:2CM}
    Let $\Delta$ be a pure simplicial complex of positive dimension and let $U$ 
    be a subset of the vertex set of $\Delta$ which has the following property: 
    every facet of $\Delta$ contains exactly one element of $U$. 
     \begin{itemize}
        \item[(i)] If $\Delta$ is Cohen-Macaulay over $\kk$ and $\link_\Delta (u)$ 
        is doubly Cohen-Macaulay over $\kk$ for every $u \in U$, then $\Delta \sm U$ 
        is doubly Cohen-Macaulay over $\kk$.  
        \item[(ii)] If $\Delta$ is Buchsbaum over $\kk$ and $\link_\Delta (u)$ is 
        doubly Cohen-Macaulay over $\kk$ for every $u \in U$, then $\Delta \sm U$ 
        is Buchsbaum* over $\kk$.  
      \end{itemize}
  \end{theorem}   

  \noindent
  \emph{Proof.}
    By \cite[Proposition 2.5 (i)]{AW}, it suffices to prove part (ii). Since $\Delta$ 
    is pure, our assumption on $U$ implies that $\Delta \sm U$ is pure as well. It also
    implies that $\Delta$ is balanced of type $(d-1, 1)$, in the sense of \cite{Sta79}, 
    where $d-1 = \dim \Delta$. Using the equivalence ${\rm (i)} \Leftrightarrow {\rm 
    (ii)}$ in Proposition \ref{prop:buchchar} and the rank-selection theorem \cite[Theorem 
    4.3]{Sta79} for balanced Cohen-Macaulay simplicial complexes, it follows as in the 
    beginning of \cite[Section 3]{BK} that $\Delta \sm U$ is Buchsbaum over $\kk$. The 
    proof of \cite[Theorem 3.1]{BK} shows (by induction on the cardinality of $U$) that 
    $\Delta \sm U$ satisfies the condition of Definition \ref{def:Buch*}, using exactly 
    the hypotheses that $\Delta$ is Buchsbaum over $\kk$ and that $\link_\Delta (u)$ is 
    doubly Cohen-Macaulay over $\kk$ for every $u \in U$. Thus the proof follows from 
    that of \cite[Theorem 3.1]{BK}. 
  \qed

  \medskip
  \noindent 
  \emph{Proof of Theorem \ref{thm:2CMposets}.} The statement follows by applying the 
  appropriate part of Theorem \ref{thm:2CM} to the order complex $\Delta (\bar{P})$ 
  and choosing $U$ as the set of maximal elements of $P$. \qed

  \section{Proof of Theorem \ref{thm:main}}
  \label{sec:proof}

  Throughout this section, $P$ is a (finite) lower Eulerian Cohen-Macaulay (over $\kk$) 
  poset of rank $d$. When $d \ge 2$, we denote by $Q$ the poset which is obtained from
  $P$ by removing all maximal elements and by $R$ the poset which is obtained from $Q$
  by removing all atoms. Thus $Q$ is lower Eulerian and graded of rank $d-1$ and $R$
  is graded of rank $d-2$. Moreover, by the rank-selection theorem for balanced 
  Cohen-Macaulay complexes \cite[Theorem 4.5]{StaCCA}, both $Q$ and $R$ are 
  Cohen-Macaulay over $\kk$.

  \begin{lemma} \label{lem:proof}
    Under the assumptions and notation of Theorem \ref{thm:main}, and assuming that 
    $P$ has rank $d \ge 2$, we have  

      \begin{eqnarray*}
      \sum_{x \in \aA(P)} |\mu_{\hat{P}} (x, \hat{1})| \, - \, \alpha_P \, |\mu_{\hat{P}} 
      (\hat{0}, \hat{1})| 
      &=& (\alpha_P - 1) \, |\widetilde{\chi} (Q)| \, - \, |\widetilde{\chi} (R)| \ + 
      \sum_{y \in \mM(P)} (\alpha(y) - \alpha_P), \\
      \end{eqnarray*}
    where $\mM (P)$ is the set of maximal elements of $P$ and $\alpha(y)$ is the number of
    atoms $x \in \aA(P)$ satisfying $x \le_P y$, for $y \in \mM(P)$.
  \end{lemma}   

  \noindent
  \emph{Proof.}
    Since $P$ is Cohen-Macaulay over $\kk$, so is $\hat{P} = P \cup \{\hat{1}\}$ and hence 
    \cite[Proposition 3.8.11]{StaEC1} we have $(-1)^d \mu_{\hat{P}} (x, \hat{1}) \ge 0$ for 
    every $x \in \aA(P)$. Using this fact and Lemma \ref{lem:chi}, we find that  
      \begin{eqnarray*}
      \sum_{x \in \aA(P)} |\mu_{\hat{P}} (x, \hat{1})| &=& (-1)^d \sum_{x \in \aA(P)} 
      \mu_{\hat{P}} (x, \hat{1}) \ = \ (-1)^d \sum_{x \in \aA(P)} \, \sum_{x \, \le_P \, y} 
      (-1)^{\rho(y)} \\
      & & \\
      &=& (-1)^d \sum_{x \in \aA(Q)} \, \sum_{x \, \le_Q \, y} (-1)^{\rho(y)} \ +  
      \sum_{y \in \mM(P)} \alpha(y), 
      \end{eqnarray*}

    \smallskip
    \noindent
    where $\rho(y)$ is the rank of $y$ in $P$. Similarly, we have
      \begin{eqnarray*}
      |\mu_{\hat{P}} (\hat{0}, \hat{1})| &=& (-1)^{d-1} \, \widetilde{\chi} (P) \ = \ 
      (-1)^{d-1} \, \psi (P) \ = \ (-1)^{d-1} \, \psi (Q) \, + \, f_{d-1} (P) \\
      &=& (-1)^{d-1} \, \widetilde{\chi} (Q) \, + \, f_{d-1} (P) \ = \ f_{d-1} (P) \, - \, 
      |\widetilde{\chi} (Q)|
      \end{eqnarray*}
    and hence

      \begin{eqnarray*}
      \sum_{x \in \aA(P)} |\mu_{\hat{P}} (x, \hat{1})| \, - \, \alpha_P \, |\mu_{\hat{P}} 
      (\hat{0}, \hat{1})| 
      &=& (-1)^d \sum_{x \in \aA(Q)} \, \sum_{x \, \le_Q \, y} (-1)^{\rho(y)} \, + \, 
          \alpha_P \, |\widetilde{\chi} (Q)| \, + \\
      & & \\
      & & \sum_{y \in \mM(P)} (\alpha(y) - \alpha_P). 
      \end{eqnarray*}
    Thus, it suffices to show that
      $$ (-1)^{d-1} \sum_{x \in \aA(Q)} \, \sum_{x \, \le_Q \, y} (-1)^{\rho(y)} \ = \ 
        |\widetilde{\chi} (Q)| \, + \, |\widetilde{\chi} (R)| $$
    or, equivalently, that 
      \begin{equation} \label{eq:step}
        \sum_{x \in \aA(Q)} \, \sum_{x \, \le_Q \, y} (-1)^{\rho(y) - 1} \ = \
        \widetilde{\chi} (Q) \, - \, \widetilde{\chi} (R). 
      \end{equation} 
    Indeed, we have
       \begin{eqnarray*}
      \sum_{x \in \aA(Q)} \, \sum_{x \, \le_Q \, y} (-1)^{\rho(y) - 1} 
      &=& \sum_{\hat{0} \ne x \, \le_Q \, y} (-1)^{\rho(y) - \rho(x)} \ - 
         \sum_{\hat{0} \ne x \, \le_R \, y} (-1)^{\rho(y) - \rho(x)} \\
      & & \\
      &=& \psi (\inte(\bar{Q})_\circ) \, - \, \psi (\inte(\bar{R})_\circ), 
      \end{eqnarray*}
    where $\inte(\bar{Q})_\circ$ (respectively, $\inte(\bar{R})_\circ$) is the poset 
    obtained from $\inte(\bar{Q})$ (respectively, $\inte(\bar{R})$) by adding a minimum 
    element. Since $\bar{Q}$ and $\bar{R}$ are locally Eulerian, the posets 
    $\inte(\bar{Q})_\circ$ and $\inte(\bar{R})_\circ$ are lower Eulerian by Proposition 
    \ref{prop:inteLE}. Thus, Lemma \ref{lem:chi} implies that $\psi(\inte(\bar{Q})_\circ) 
    = \widetilde{\chi} (\inte(\bar{Q})_\circ)$ and $\psi(\inte(\bar{R})_\circ) = 
    \widetilde{\chi} (\inte(\bar{R})_\circ)$. Finally, we note that $\widetilde{\chi} 
    (\inte(\bar{Q})_\circ) = \widetilde{\chi} (Q)$ and $\widetilde{\chi} 
    (\inte(\bar{R})_\circ) = \widetilde{\chi} (R)$ by Proposition \ref{prop:inte} and 
    hence (\ref{eq:step}) follows. 
  \qed

  \medskip
  We now proceed with the proof of Theorem \ref{thm:main}. We recall that for every 
  atom $x$ of $Q$ we have the natural map $\rho_x: \widetilde{H}_{d-2} (\Delta(\bar{Q}); 
  \kk) \rightarrow \widetilde{H}_{d-3} (\Delta(Q_{> x}); \kk)$, shown to be surjective
  in Corollary \ref{cor:surjectivity}.

  \medskip
  \noindent 
  \emph{Proof of Theorem \ref{thm:main}.} Let $d$ be the rank of $P$. We assume that 
  $d \ge 2$, since the result is trivial otherwise. By Corollary \ref{cor:Qhomology}, 
  we may choose $\bB(P) \subseteq \mM(P)$ so that the classes $\omega(y)$ for $y \in 
  \bB(P)$ form a basis of the $\kk$-vector space $\widetilde{H}_{d-2} (\Delta(\bar{Q}); 
  \kk)$. Since $Q$ is Cohen-Macaulay over $\kk$, it follows from (\ref{eq:chiCM}) 
  that the cardinality of $\bB(P)$ is equal to $|\widetilde{\chi} (Q)|$. Hence 
  Lemma \ref{lem:proof} implies that  

    \begin{eqnarray*}
    \sum_{x \in \aA(P)} |\mu_{\hat{P}} (x, \hat{1})| \, - \, \alpha_P \, |\mu_{\hat{P}} 
    (\hat{0}, \hat{1})| 
    & \ge & (\alpha_P - 1) \, |\widetilde{\chi} (Q)| \, - \, |\widetilde{\chi} (R)| \ + 
    \sum_{y \in \bB(P)} (\alpha(y) - \alpha_P) \\
    &=& \sum_{y \in \bB(P)} (\alpha(y) - 1) \, - \, |\widetilde{\chi} (R)|.
    \end{eqnarray*}
  Therefore, it suffices to show that 
    \begin{equation} \label{eq:suff}
      |\widetilde{\chi} (R)| \ \le \ \sum_{y \in \bB(P)} (\alpha(y) - 1).
    \end{equation} 
  For $x \in \aA(Q)$ we consider the natural map $\rho_x: \widetilde{H}_{d-2} 
  (\Delta(\bar{Q}); \kk) \rightarrow \widetilde{H}_{d-3} (\Delta(Q_{> x}); \kk)$ and the 
  map $\widetilde{H}_{d-3} (\Delta(Q_{> x}); \kk) \rightarrow \widetilde{H}_{d-3} 
  (\Delta(\bar{R}); \kk)$, induced by inclusion. These maps yield the sequence of linear 
  maps
    \begin{equation} \label{eq:seq}
      \bigoplus_{x \in \aA(Q)} \widetilde{H}_{d-2} (\Delta(\bar{Q}); \kk) \ \rightarrow  
       \bigoplus_{x \in \aA(Q)}  \widetilde{H}_{d-3} (\Delta(Q_{> x}); \kk) \ \rightarrow   
       \ \widetilde{H}_{d-3} (\Delta(\bar{R}); \kk).
    \end{equation} 
  We observe that both maps in this sequence are surjective, the one on the left since
  every map $\rho_x$ is surjective by Corollary \ref{cor:surjectivity}, and the one on 
  the right by Lemma \ref{lem:Qhomology} (applied to the dual of $\bar{Q}$). For $y \in 
  \bB(P)$, let us denote by $P(y)$ the poset which is obtained from the open interval 
  $(\hat{0}, y)$ of $P$ by removing all its minimal elements. Since there are exactly 
  $\alpha(y)$ such elements and $P(y)$ is Cohen-Macaulay of rank $d-3$, we have 
    $$ \dim_\kk \widetilde{H}_{d-3} (\Delta(P(y)); \kk) \ = \ \alpha(y) - 1 $$
  by (\ref{eq:chiCM}) and Lemma \ref{lem:mobiusminusatoms}. Clearly, for $x \in \aA(P)$ 
  and $y \in \bB(P)$ we may have $\rho_x (\omega(y)) \ne 0$ only if $x \le_P y$. This 
  implies that the image of the composition of the two maps in (\ref{eq:seq}) is 
  contained in (the image of the map induced by the inclusion of) the sum of the spaces 
  $\widetilde{H}_{d-3} (\Delta(P(y)); \kk)$ for $y \in \bB(P)$. This fact and the 
  surjectivity of the maps in (\ref{eq:seq}) imply that 
    $$ |\widetilde{\chi} (R)| \ = \ \dim_\kk \widetilde{H}_{d-3} (\Delta(\bar{R}); \kk) 
       \ \le \ \sum_{y \in \bB(P)} \dim_\kk \widetilde{H}_{d-3} (\Delta(P(y)); \kk) \ = 
       \ \sum_{y \in \bB(P)} (\alpha(y) - 1). $$
  This shows (\ref{eq:suff}) and completes the proof of the theorem.
  \qed 

  \section{Applications}
  \label{sec:app}

  This section deduces Corollaries \ref{cor:cubical} and \ref{cor:toric} from Theorem
  \ref{thm:main}. 

  \subsection{The toric $h$-vector}
  \label{sec:toric}

  We first recall from \cite{Sta87} the definition of the toric $h$-vector of a (finite)
  lower Eulerian poset $P$. We denote by $\hat{0}$ the minimum element of $P$, as usual, 
  and by $d$ the maximum rank $\rho(y)$ of an element $y \in P$. For $y \in P$ we define 
  two polynomials $f(P, y; q)$ and $g(P, y; q)$ by the following rules: 

  \smallskip
    \begin{itemize}
    \itemsep=2pt
      \item[(a)] $f(P, \hat{0}; q) = g(P, \hat{0}; q) = 1$. 
      \item[(b)] If $y \in \bar{P}$ and $f(P, y; q) = k_0 + k_1 q + \cdots$, then
        \begin{equation} \label{eq:deftoric-g}
          g(P, y; q) \ = \ k_0 + \sum_{i=1}^m \ (k_i - k_{i-1}) q^i, 
        \end{equation}
      where $m = \lfloor (\rho(y) - 1) / 2 \rfloor$. 
      \item[(c)] If $z \in \bar{P}$, then 
        \begin{equation} \label{eq:deftoric-f}
          f(P, z; q) \ = \ \sum_{y <_P z} \ g(P, y; q) (q-1)^{\rho(y, z) - 1}.
        \end{equation}
    \end{itemize}

  \noindent
  The toric $h$-vector of $P$ is the sequence $h(P) = (h_0 (P), h_1 (P),\dots,h_d 
  (P))$ defined by
    \begin{equation} \label{eq:deftoric}
      h_d (P) + h_{d-1} (P) q + \cdots + h_0 (P) q^d \ = \ \sum_{y \in P} \ g(P, y; q) 
      (q-1)^{d - \rho(y)}. 
    \end{equation}
  For instance, we have $h_0 (P) = 1$, $h_1 (P) = f_0 (P) - d$ and $h_d (P) = (-1)^{d-1} 
  \widetilde{\chi} (P)$, where $f_0 (P)$ is the number of atoms of $P$ (see \cite[p. 
  198]{Sta87}). These formulas imply that $h_1 (P) \ge 0$, if $P$ is a meet-semilattice 
  (meaning that any two elements of $P$ have a greatest lower bound) and $h_d (P) \ge 0$, 
  if $P$ is Cohen-Macaulay over $\kk$. We will write $h(P) \ge 0$ if we have $h_i (P) \ge 
  0$ for $0 \le i \le d$. The following conjecture is part of \cite[Conjecture 4.2]{Sta87}. 

  \begin{conjecture} {\rm (\cite[Conjecture 4.2 (b)]{Sta87})} \label{conj:toric}
    For every lower Eulerian Cohen-Macaulay meet-semilattice $P$ we have $h(P) \ge 0$. 
  \end{conjecture}

  The generalized Dehn-Somerville equations \cite[Theorem 2.4]{Sta87} assert that the 
  polynomial $f(P, z; q)$ is symmetric of degree $\rho(z) - 1$ for every $z \in \bar{P}$.
  Thus, writing $f(P, z; q) = \sum_{i=0}^r k_i q^i$ with $r = \rho(z) - 1$, we have $k_i 
  = k_{r - i}$ for $0 \le i \le r$. 

  We now deduce from Theorem \ref{thm:main} that $h_{d-1} (P) \ge 0$ for every lower 
  Eulerian Cohen-Macaulay meet-semilattice $P$ of rank $d$.  

  \medskip
  \noindent 
  \emph{Proof of Corollary \ref{cor:toric}.} We note that for every maximal element 
  $y \in P$, the interval $[\hat{0}, y]$ in $P$ is an Eulerian lattice of rank $d$ and 
  therefore has at least $d$ atoms (this holds, more generally, for all graded lattices 
  of rank $d$ with nowhere vanishing M\"obius function). Thus, in the notation of 
  Theorem \ref{thm:main}, we have $\alpha_P \ge d$. We claim that
    \begin{equation} \label{eq:htoric{d-1}}
      h_{d-1} (P) \ = \ \sum_{x \in \aA(P)} \, |\widetilde{\chi} (P_{\ge x})| \, - \, 
      d \ |\widetilde{\chi} (P)|,
    \end{equation}
  where $P_{\ge x} := \{y \in P: x \le_P y\}$ is considered as a subposet of $P$. In
  view of (\ref{eq:main}) and the inequality $\alpha_P \ge d$, it suffices to prove
  the claim. For $y \in P$, let $\alpha(y)$ denote the number of atoms of $P$ in the 
  interval $[\hat{0}, y]$. The constant term of $g(P, y; q)$ is equal to 1 for every 
  $y \in P$. Furthermore, it follows from equations (\ref{eq:deftoric-g}) and 
  (\ref{eq:deftoric-f}) and from the symmerty of the polynomials $f(P, z; q)$ that 
  the coefficient of $q$ in $g(P, y; q)$ is equal to $\alpha(y) - \rho(y)$. Thus, 
  equation (\ref{eq:deftoric}) implies that
    $$ h_{d-1} (P) \ = \ \sum_{y \in P} \, (-1)^{d - \rho(y)} \, (\alpha(y) - d). $$
  An easy computation, essentially already carried out in the proof of Lemma 
  \ref{lem:proof}, shows that the right-hand side of the previous equation is equal
  to the right-hand side of (\ref{eq:htoric{d-1}}). This completes the proof.
  \qed

  \subsection{The cubical $h$-vector}
  \label{sec:cubical}

  Cubical posets form an important class of lower Eulerian posets. As mentioned in the 
  introduction, a cubical poset is a (finite) poset having a minimum element $\hat{0}$, 
  such that for every $x \in P$, the interval $[\hat{0}, x]$ in $P$ is isomorphic to 
  the poset of faces of a cube (the dimension of which has to equal one less than the 
  rank of $x$). We assume that $P$ is graded of rank $d$ and denote by $f_{i-1} (P)$ 
  the number of elements of $P$ of rank $i$ for $0 \le i \le d$, as usual. The 
  \emph{cubical $h$-polynomial} of $P$ can be defined (see \cite{Ad96}) by the formula 
    \begin{equation} \label{eq:cubicalp}
      (1+q) h^{(c)} (P, q) \ = \ 2^{d-1} \, + \, q \, h^{(sc)} (P, q) \, + \, 
      (-2)^{d-1} \, \widetilde{\chi} (P) \, q^{d+1},
    \end{equation}
  where   
    \begin{equation} \label{eq:scubicalp}
      h^{(sc)} (P, q) \ = \ \sum_{i=0}^{d-1} \, f_i (P) \, (2q)^i (1-q)^{d-i-1} 
    \end{equation}
  is the short cubical $h$-polynomial of $P$ and $\widetilde{\chi} (P) = \psi(P) = 
  \mu_{\hat{P}} (\hat{0}, \hat{1})$ (see Section \ref{subsec:lower}). The function $h^{(c)} 
  (P, q)$ is a polynomial in $q$ of degree at most $d$. The \emph{cubical $h$-vector} of 
  $P$ is defined as the sequence $(h^{(c)}_0 (P), h^{(c)}_1 (P),\dots,h^{(c)}_d (P))$, 
  where 
    $$ h^{(c)} (P, q) \ = \ \sum_{i=0}^d h^{(c)}_i (P) q^i. $$
  Clearly, (\ref{eq:cubicalp}) implies that $h^{(c)}_d (P) = (-2)^{d-1} \widetilde{\chi} 
  (P)$ and hence we have $h^{(c)}_d (P) \ge 0$, if $P$ is Cohen-Macaulay over $\kk$. By 
  direct computation we also find that 
    \begin{equation} \label{eq:h_{d-1}}
      h^{(c)}_{d-1} (P) \ = \ (-2)^{d-1} \, + \ \sum_{i=1}^d \ (-1)^{d-i-1}
      \, (2^{d-1} - 2^{i-1}) f_{i-1} (P).   
    \end{equation}
  For instance, we have 
    $$ h^{(c)}_{d-1} (P) \ = \ \cases{ f_0 (P) - 2, & if \ $d=2$ \cr
       2 f_1 (P) - 3 f_0 (P) + 4, & if \ $d=3$ \cr
       4 f_2 (P) - 6 f_1 (P) + 7 f_0 (P) - 8, & if \ $d=4$. } $$
  We now deduce from Theorem \ref{thm:main} that $h^{(c)}_{d-1} (P) \ge 0$ for every 
  Cohen-Macaulay cubical poset $P$ of rank $d$.  

  \medskip
  \noindent 
  \emph{Proof of Corollary \ref{cor:cubical}.} Since the poset $P$ is cubical (and 
  graded) of rank $d$, in the notation of Theorem \ref{thm:main} we have $\alpha_P = 
  2^{d-1}$. Thus, in view of (\ref{eq:main}) and since $P$ is Cohen-Macaulay over $\kk$, 
  it suffices to verify (\ref{eq:mainc}). Comparing the coefficients of $x^d$ in the 
  two sides of (\ref{eq:cubicalp}) we get 
    \begin{equation} \label{eq:h_{d-1}step}
      h^{(c)}_{d-1} (P) \ = \ h^{(sc)}_{d-1} (P) - h^{(c)}_d (P) \ = \ h^{(sc)}_{d-1} 
      (P) \, - \, (-2)^{d-1} \, \widetilde{\chi} (P), 
    \end{equation}
  where $h^{(sc)}_{d-1} (P)$ denotes the coefficient of $x^{d-1}$ in $h^{(sc)} (P, q)$. 
  We next recall that for $x \in \bar{P}$, the subposet $P_{\ge x} := \{y \in P: x 
  \le_P y\}$ of $P$ is a simplicial poset and that, by a fundamental observation of 
  Hetyei (see \cite[Theorem 9]{Ad96}), we have 
    \begin{equation} \label{eq:scubicalh}
      h^{(sc)} (P, q) \ = \ \sum_{x \in \aA(P)} h(P_{\ge x}, q),
    \end{equation}
  where $h(Q, q)$ stands for the simplicial $h$-vector of a simplicial poset $Q$ (as 
  defined, for instance, by the right-hand side of (\ref{eq:h})). Equation 
  (\ref{eq:scubicalh}) implies that
    \begin{equation} \label{eq:hchiP}
      h^{(sc)}_{d-1} (P) \ = \ \sum_{x \in \aA(P)} (-1)^d \, \widetilde{\chi} (P_{\ge x}). 
    \end{equation}
  Finally, (\ref{eq:h_{d-1}step}) and (\ref{eq:hchiP}) imply (\ref{eq:mainc}) and the 
  proof follows. \qed

  \medskip
  Adin \cite[Question 1]{Ad96} raised the question whether $h^{(c)} (P, q) \ge 0$ holds 
  for every Cohen-Macaulay cubical meet-semilattice $P$ (where the inequality is meant to 
  hold coefficientwise). In view of Corollary \ref{cor:cubical}, it is natural to extend 
  this question as follows.  

  \begin{question} \label{que:cubical}
    Does $h^{(c)} (P, q) \ge 0$ hold for every Cohen-Macaulay cubical poset $P$? 
  \end{question}

  \section*{Acknowledgments} The author wishes to thank the anonymous referees for a 
  careful reading of the manuscript and for a number of useful suggestions on the 
  presentation.

  \end{document}